\input amstex
\documentstyle{amsppt}
%
\catcode`@=11
\redefine\output@{%
  \def\break{\penalty-\@M}\let\par\endgraf
  \ifodd\pageno\global\hoffset=105pt\else\global\hoffset=8pt\fi  
  \shipout\vbox{%
    \ifplain@
      \let\makeheadline\relax \let\makefootline\relax
    \else
      \iffirstpage@ \global\firstpage@false
        \let\rightheadline\frheadline
        \let\leftheadline\flheadline
      \else
        \ifrunheads@ 
        \else \let\makeheadline\relax
        \fi
      \fi
    \fi
    \makeheadline \pagebody \makefootline}%
  \advancepageno \ifnum\outputpenalty>-\@MM\else\dosupereject\fi
}
\def\Beta{\mathchar"0\hexnumber@\rmfam 42}
\catcode`\@=\active
\nopagenumbers
\chardef\textvolna='176

\chardef\bigalpha='013
\def\negskp{\hskip -2pt}

\chardef\degree="5E
\def\compos{\,\raise 1pt\hbox{$\sssize\circ$} \,}

\def\blue#1{#1}

\catcode`#=11\def\diez{#}\catcode`#=6
\catcode`&=11\catcode`&=4
\catcode`_=11\def\podcherkivanie{_}\catcode`_=8
\catcode`\^=11\catcode`\^=7
\catcode`~=11\def\volna{~}\catcode`~=\active
\def\mycite#1{\cite{\blue{#1}}\immediate\special{ps:
     ShrHPSdict begin /ShrBORDERthickness 0 def}}
\def\myciterange#1#2#3#4{\cite{\blue{#2#3#4}}\immediate\special{ps:
     ShrHPSdict begin /ShrBORDERthickness 0 def}}
\def\mytag#1{%
    \tag#1}
\def\mythetag#1{\thetag{\blue{#1}}\immediate\special{ps:
     ShrHPSdict begin /ShrBORDERthickness 0 def}}
\def\myrefno#1{\no#1}
\def\myhref#1#2{\blue{#2}\immediate\special{ps:
     ShrHPSdict begin /ShrBORDERthickness 0 def}}
\def\myEarXivlink{\myhref{http://arXiv.org}{http:/\negskp/arXiv.org}}

\def\mytheorem#1{\csname proclaim\endcsname{Theorem #1}}
\def\mytheoremwithtitle#1#2{\csname proclaim\endcsname{Theorem #1#2}}
\def\mythetheorem#1{\blue{#1}\immediate\special{ps:
     ShrHPSdict begin /ShrBORDERthickness 0 def}}
\def\mylemma#1{\csname proclaim\endcsname{Lemma #1}}
\def\mylemmawithtitle#1#2{\csname proclaim\endcsname{Lemma #1#2}}
\def\mythelemma#1{\blue{#1}\immediate\special{ps:
     ShrHPSdict begin /ShrBORDERthickness 0 def}}
\def\mycorollary#1{\csname proclaim\endcsname{Corollary #1}}

\def\mydefinition#1{\definition{Definition #1}}
\def\mythedefinition#1{\blue{#1}\immediate\special{ps:
     ShrHPSdict begin /ShrBORDERthickness 0 def}}
\def\myconjecture#1{\csname proclaim\endcsname{Conjecture #1}}
\def\myconjecturewithtitle#1#2{\csname proclaim\endcsname{Conjecture #1#2}}

\def\myproblem#1{\csname proclaim\endcsname{Problem #1}}
\def\myproblemwithtitle#1#2{\csname proclaim\endcsname{Problem #1#2}}


\pagewidth{360pt}
\pageheight{606pt}
\topmatter
\title
A note on invertible quadratic transformations 
of the real plane.
\endtitle
\rightheadtext{A note on invertible quadratic transformations \dots}
\author
Ruslan Sharipov
\endauthor
\address Bashkir State University, 32 Zaki Validi street, 450074 Ufa, Russia
\endaddress
\email
\myhref{mailto:r-sharipov\@mail.ru}{r-sharipov\@mail.ru}
\endemail
\abstract
     A polynomial transformation of the real plane $\Bbb R^2$ is a mapping
$\Bbb R^2\to\Bbb R^2$ given by two polynomials of two variables. Such a
transformation is called quadratic if the degrees of its polynomials are not
greater than two. In the present paper an exhaustive description of invertible 
quadratic transformations of the real plane is given. Their application to
the perfect cuboid problem is discussed. 
\endabstract
\endtopmatter
\TagsOnRight
\document

\head
1. Introduction.
\endhead
     A perfect cuboid is a rectangular parallelepiped whose edges, whose face 
diagonals and whose space diagonal all are of integer lengths. None of such 
cuboids is found thus far. The problem of finding them or proving their 
non-existence is still open, see its history in \myciterange{1}{1}{--}{47}.
\par
     Slanted (non-rectangular) perfect cuboids are known. Some of them were
found in \mycite{48}. Infinite families of slanted perfect cuboids were found 
in \mycite{49} and \mycite{50}.\par
     In some approaches the search for rectangular perfect cuboids is reduced 
to a single Diophantine equation. One of such Diophantine equations was derived 
in \mycite{51}. On the basis of this equation in \mycite{52} three cuboid 
conjectures were formulated. The first cuboid conjecture is rather simple. 
Though it was not yet proved, in \mycite{53} it was shown that there no 
perfect cuboids in the case of the first cuboid conjecture.\par
     The second and the third cuboid conjectures are more complicated. They 
were considered in \mycite{54} and \mycite{55}, but were not proved. At the 
present time the case of the second cuboid conjecture is being intensively 
studied using the asymptotic approach (see \myciterange{56}{56}{--}{59}). The 
approach used in \myciterange{60}{60}{--}{72} is quite different. We do not 
consider this approach below in the present paper.\par
     In the case of the second cuboid conjecture the search for rectangular 
perfect cuboids is reduced to the following Diophantine equation of
tenth degree with respect to the positive integer variable $t>0$: 
$$
\gathered
t^{10}+(2\,q^{\kern 0.7pt 2}+p^{\kern 1pt 2})\,(3
\,q^{\kern 0.7pt 2}-2\,p^{\kern 1pt 2})\,t^8+(q^{\kern 0.5pt 8}+10
\,p^{\kern 1pt 2}\,q^{\kern 0.5pt 6}+\\
+\,4\,p^{\kern 1pt 4}\,q^{\kern 0.5pt 4}-14\,p^{\kern 1pt 6}\,q^{\kern 0.7pt 2}
+p^{\kern 1pt 8})\,t^6-p^{\kern 1pt 2}\,q^{\kern 0.7pt 2}
\,(q^{\kern 0.5pt 8}-14\,p^{\kern 1pt 2}\,q^{\kern 0.5pt 6}
+4\,p^{\kern 1pt 4}\,q^{\kern 0.5pt 4}+\\
+\,10\,p^{\kern 1pt 6}\,q^{\kern 0.7pt 2}+p^{\kern 1pt 8})\,t^4
-p^{\kern 1pt 6}\,q^{\kern 0.5pt 6}\,(q^{\kern 0.7pt 2}
+2\,p^{\kern 1pt 2})\,(3\,p^{\kern 1pt 2}-2\,q^{\kern 0.7pt 2})\,t^2
-q^{\kern 0.7pt 10}\,p^{\kern 1pt 10}=0.
\endgathered\quad
\mytag{1.1}
$$
Two coprime positive integer numbers $p\neq q$ are parameters of the equation
\mythetag{1.1}. They define a point $(p,q)$ on the $p\,q$\,-\,coordinate
plane $\Bbb R^2$.\par
    Within the asymptotic approach to the equation 
\mythetag{1.1} its parameters $p$ and $q$ tend to infinity either separately
like in \mycite{56} and \mycite{57} or simultaneously like in \mycite{58} and 
\mycite{59}. In order to arrange a proper passage to the limit $(p,q)\to\infty$ 
in \mycite{58} and \mycite{59} linear and nonlinear transformations of the form
$$
\xalignat 2
&\hskip -2em
\tilde p=\tilde p(p,q),
&&\tilde q=\tilde q(p,q).
\mytag{1.2}
\endxalignat
$$
were used. There are two natural requirements for the transformations 
\mythetag{1.2} applied to the parameters of the Diophantine equation 
\mythetag{1.1}: 
\roster
\item"1)" they should be invertible;
\item"2)" they should map bijectively the integer $p\,q$\,-\,greed onto
itself. 
\endroster
Polynomial transformations are the best candidates for this role. For example
in \mycite{59} the following cubic transformations of the form \mythetag{1.2}
were used:
$$
\xalignat 2
&\hskip -2em
\tilde p=B\,q^{\kern 0.7pt 3}-p,
&&\tilde q=q.
\mytag{1.3}
\endxalignat
$$
Unfortunately the coefficient $B$ in \mythetag{1.3} is restricted to 
integer numbers. Therefore the transformations \mythetag{1.3} appeared 
to be insufficient for to cover all needs they were designed to serve
in \mycite{59}.\par
     The main goal of the present paper is to extend the number of
invertible polynomial transformations applicable to Diophantine equation 
\mythetag{1.1}. We consider quadratic transformations which are most
simple after linear ones.\par
     It is worth to note that complex polynomial transformations in $\Bbb C^n$
are associated with the Jacobian conjecture, which is another open problem
in mathematics (see \mycite{73}). The corresponding real Jacobian conjecture
in $\Bbb R^n$ is invalid. A counter example to it was given by S.~I.~Pinchuk
in \mycite{74}.\par
\head
2. Quadratic transformations.
\endhead
     A general quadratic transformation of $\Bbb R^2$ is given by the following
formulas:
$$
\hskip -3em
\aligned
&\tilde p=a_{20}\,p^2+2\,a_{11}\,p\,q+a_{02}\,q^2+2\,a_{10}\,p
+2\,a_{01}\,q+a_{00},\\
&\tilde q=b_{20}\,p^2+2\,b_{11}\,p\,q+b_{02}\,q^2+2\,b_{10}\,p
+2\,b_{01}\,q+b_{00}.
\endaligned
\mytag{2.1}
$$
Without loss of generality one can assume that 
$$
\hskip -2em
b_{11}=0.
\mytag{2.2}
$$
Indeed, if the condition \mythetag{2.2} is not fulfilled and if $a_{11}\neq 0$, 
we can compose \mythetag{2.1} with the following invertible linear transformation:
$$
\xalignat 2
&\hskip -2em
\hat p=\tilde p,
&&\hat q=\tilde q-\frac{b_{11}}{a_{11}}\,\tilde p.
\mytag{2.3}
\endxalignat
$$
If $a_{11}=0$, we can just exchange variables by applying the transformation
$$
\xalignat 2
&\hskip -2em
\hat p=\tilde q,
&&\hat q=\tilde p
\mytag{2.4}
\endxalignat
$$
instead of \mythetag{2.3}. Like \mythetag{2.3}, the transformation \mythetag{2.4}
is invertible. In both cases the composite transformation will obey the condition
\mythetag{2.2}.\par
      Assume that the coefficients $b_{20}$ and $b_{11}$ are proportional to the
coefficients $a_{20}$ and $a_{11}$ in \mythetag{2.1}, i\.\,e\. assume that they
obey the relationship
$$
\hskip -2em
\frac{b_{20}}{b_{11}}=\frac{a_{20}}{a_{11}}.
\mytag{2.5}
$$
In this case, composing \mythetag{2.3} with \mythetag{2.1}, we get a quadratic
transformation of the form \mythetag{2.1} with two coefficients $b_{20}$ and 
$b_{11}$ being equal to zero:
$$
\xalignat 2
&\hskip -2em
b_{20}=0,
&&b_{11}=0.
\mytag{2.6}
\endxalignat
$$
Similarly, if the relationship 
$$
\hskip -2em
\frac{b_{02}}{b_{11}}=\frac{a_{02}}{a_{11}}.
\mytag{2.7}
$$
is fulfilled, then upon composing \mythetag{2.3} with \mythetag{2.1} we shall
have 
$$
\xalignat 2
&\hskip -2em
b_{11}=0,
&&b_{02}=0.
\mytag{2.8}
\endxalignat
$$
Apart from \mythetag{2.5} and \mythetag{2.7}, there is a third option:
$$
\hskip -2em
\frac{b_{02}}{b_{20}}=\frac{a_{02}}{a_{20}}.
\mytag{2.9}
$$
In this case instead of \mythetag{2.3} we choose the following
linear transformation:  
$$
\xalignat 2
&\hskip -2em
\hat p=\tilde p,
&&\hat q=\tilde q-\frac{b_{20}}{a_{20}}\,\tilde p.
\mytag{2.10}
\endxalignat
$$
Provided \mythetag{2.9} is fulfilled, upon composing \mythetag{2.10} with 
\mythetag{2.1} we shall have 
$$
\xalignat 2
&\hskip -2em
b_{20}=0,
&&b_{02}=0.
\mytag{2.11}
\endxalignat
$$
The relationships \mythetag{2.5}, \mythetag{2.7}, and \mythetag{2.9} can 
be written in a denominator-free form:
$$
\xalignat 3
&\det\Vmatrix a_{20} & a_{11}\!\\
\vspace{1ex}
b_{20} & b_{11}\!\endVmatrix=0,
&&\det\Vmatrix a_{11} & a_{02}\!\\
\vspace{1ex}
b_{11} & b_{02}\!\endVmatrix=0,
&&\det\Vmatrix a_{20} & a_{02}\!\\
\vspace{1ex}
b_{20} & b_{02}\!\endVmatrix=0.
\qquad\quad
\mytag{2.12}
\endxalignat
$$\par
      Typically the relationships \mythetag{2.12} are not fulfilled. However,
we can use a linear transformation in order to change the coefficients 
in \mythetag{2.1}. Let's set
$$
\xalignat 2
&\hskip -2em
p=c_{11}\,\check p+c_{12}\,\check q,
&&q=c_{21}\,\check p+c_{22}\,\check q.
\mytag{2.13}
\endxalignat
$$
Substituting \mythetag{2.9} into \mythetag{2.1}, we get a transformation 
of the same form like \mythetag{2.1}:
$$
\aligned
&\tilde p=\check a_{20}\,\check p^2+2\,\check a_{11}\,\check p\,\check q
+\check a_{02}\,\check q^2+2\,\check a_{10}\,\check p
+2\,\check a_{01}\,\check q+\check a_{00},\\
&\tilde q=\check b_{20}\,\check p^2+2\,\check b_{11}\,\check p\,\check q
+\check b_{02}\,\check q^2+2\,\check b_{10}\,\check p
+2\,\check b_{01}\,\check q+\check b_{00}.
\endaligned
$$
Its coefficients depend on $c_{11}$, $c_{12}$, $c_{21}$, $c_{22}$ and on
the coefficients of \mythetag{2.1}. Upon substituting $\check a_{20}$,
$\check a_{11}$, $\check a_{02}$ and $\check b_{20}$, $\check b_{11}$, 
$\check b_{02}$ for $a_{20}$, $a_{11}$, $a_{02}$ and $b_{20}$, $b_{11}$, 
$b_{02}$ into \mythetag{2.8} we derive the following three relationships:
$$
\align
&\hskip -3em
\det C\,\biggl(c_{11}^2\,\vmatrix a_{20} & a_{11}\!\\
\vspace{1ex}
b_{20} & b_{11}\!\endvmatrix+c_{11}\,\vmatrix a_{20} & a_{02}\!\\
\vspace{1ex}
b_{20} & b_{02}\!\endvmatrix\,c_{21}+\vmatrix a_{11} & a_{02}\!\\
\vspace{1ex}
b_{11} & b_{02}\!\endvmatrix\,c_{21}^2\biggr)=0,\kern -2em
\mytag{2.14}\\
\vspace{2ex}
&\hskip -3em
\det C\,\biggl(c_{22}^2\,\vmatrix a_{20} & a_{11}\!\\
\vspace{1ex}
b_{20} & b_{11}\!\endvmatrix+c_{22}\,\vmatrix a_{20} & a_{02}\!\\
\vspace{1ex}
b_{20} & b_{02}\!\endvmatrix\,c_{12}+\vmatrix a_{11} & a_{02}\!\\
\vspace{1ex}
b_{11} & b_{02}\!\endvmatrix\,c_{21}^2\biggr)=0,\kern -2em
\mytag{2.15}\\
\vspace{2ex}
&\hskip -3em
\aligned
&\det C\,\biggl(c_{11}\,\vmatrix a_{20} & a_{11}\!\\
\vspace{1ex}
b_{20} & b_{11}\!\endvmatrix\,c_{12}+c_{11}\,\vmatrix a_{20} & a_{02}\!\\
\vspace{1ex}
b_{20} & b_{02}\!\endvmatrix\,c_{22}\,+\\
&\kern 7.7em+\,c_{21}\,\vmatrix a_{20} & a_{02}\!\\
\vspace{1ex}
b_{20} & b_{02}\!\endvmatrix\,c_{12}+c_{21}\,\vmatrix a_{11} & a_{02}\!\\
\vspace{1ex}
b_{11} & b_{02}\!\endvmatrix\,c_{22}\biggr)=0.
\endaligned
\mytag{2.16}
\endalign
$$
Here $\det C$ is the determinant of the following matrix:
$$
\hskip -2em
C=\Vmatrix c_{11} & c_{12}\!\\
\vspace{1ex}
 c_{21} & c_{22}\!\endVmatrix.
\mytag{2.17} 
$$
The matrix \mythetag{2.17} is non-degenerate since the transformation 
\mythetag{2.13} should be invertible. Hence $\det C\neq 0$ and we can cancel
it in \mythetag{2.14}, \mythetag{2.15}, and \mythetag{2.16}.\par
     Looking at \mythetag{2.14} and \mythetag{2.15}, we define the quadratic 
form\footnotemark
$$
\hskip -2em
\omega_1(\bold c)=c_1^2\,\vmatrix a_{20} & a_{11}\!\\
\vspace{1ex}
b_{20} & b_{11}\!\endvmatrix+c_1\,\vmatrix a_{20} & a_{02}\!\\
\vspace{1ex}
b_{20} & b_{02}\!\endvmatrix\,c_2+\vmatrix a_{11} & a_{02}\!\\
\vspace{1ex}
b_{11} & b_{02}\!\endvmatrix\,c_2^2,
\mytag{2.18} 
$$
where $\bold c\in\Bbb R^2$ is a vector with two components $c_1$ and $c_2$. Any
quadratic form in a real vector space can be definite, semi-definite, or
indefinite (see \mycite{75}).\par
\footnotetext{\ Regular quadratic forms are tensors of the type $(0,2)$. The 
quadratic form $\omega_1(\bold c)$ is a pseudo-tensor of the type $(0,2)$ in 
the sense of Definition~2.1 in \mycite{76}.} 
\head
3. The case where $\omega_1(\bold c)$ is zero.
\endhead
    In this case all of the three determinants in \mythetag{2.18} are zero. 
This means that the coefficients $b_{20}$, $b_{11}$, $b_{02}$ are proportional
to the coefficients $a_{20}$, $a_{11}$, $a_{02}$. Hence, applying some tricks
like \mythetag{2.3} and \mythetag{2.4}, we can bring \mythetag{2.1} to the 
form 
$$
\hskip -3em
\aligned
&\tilde p=a_{20}\,p^2+2\,a_{11}\,p\,q+a_{02}\,q^2+2\,a_{10}\,p
+2\,a_{01}\,q+a_{00},\\
&\tilde q=2\,b_{10}\,p+2\,b_{01}\,q+b_{00}.
\endaligned
\mytag{3.1}
$$
Since the transformation \mythetag{3.1} is assumed to be invertible, at least
one of the two coefficients $b_{10}$ and $b_{01}$ is nonzero. Therefore, applying 
some linear transformation of the form \mythetag{2.13} and a shift of origin, we can 
further simplify the formulas \mythetag{3.1}:
$$
\hskip -3em
\aligned
&\tilde p=a_{20}\,p^2+2\,a_{11}\,p\,q+a_{02}\,q^2+2\,a_{10}\,p
+2\,a_{01}\,q,\\
&\tilde q=q.
\endaligned
\mytag{3.2}
$$\par
     Assume that $a_{20}\neq 0$ and assume that $q$ is fixed. In this 
case the quadratic polynomial in \mythetag{3.2} takes some values twice 
and does not take some other values at all. This fact contradicts the 
invertibility of the transformation \mythetag{3.2}. Hence $a_{20}=0$ and 
the formulas \mythetag{3.2} simplify to the following ones:
$$
\xalignat 2
&\hskip -3em
\tilde p=2\,a_{11}\,p\,q+a_{02}\,q^2+2\,a_{10}\,p+2\,a_{01}\,q,
&&\tilde q=q.
\mytag{3.3}
\endxalignat
$$
If $a_{11}\neq 0$, then for some fixed $q$ the first polynomial in 
\mythetag{3.3} does not depend on $p$. This fact contradicts the invertibility 
of the transformation \mythetag{3.3}. Hence $a_{11}=0$. The coefficient $2\,a_{10}$ 
in \mythetag{3.3} is nonzero. Applying a linear transformation it can be reduced 
to the unity. Therefore the formula \mythetag{3.3} simplifies to 
$$
\xalignat 2
&\hskip -3em
\tilde p=p+a_{02}\,q^2+2\,a_{01}\,q,
&&\tilde q=q.
\mytag{3.4}
\endxalignat
$$
The coefficient $a_{02}$ in \mythetag{3.4} is nonzero. Applying scaling 
transformations (which are linear) in $p$, $q$, $\tilde p$, $\tilde q$
and some origin shifts, we can bring $a_{02}$ to the unity and can annul 
the coefficient $2\,a_{01}$. As a result \mythetag{3.4} turns to 
$$
\xalignat 2
&\hskip -3em
\tilde p=p+q^2,
&&\tilde q=q.
\mytag{3.5}
\endxalignat
$$
This result is formulated as the following theorem. 
\mytheorem{3.1} In the case where the associated quadratic form 	
\mythetag{2.18} is zero any invertible quadratic transformation 
\mythetag{2.1} reduces to the form \mythetag{3.5} at the expense of 
composing it with linear transformations and origin shifts. 
\endproclaim
\head
4. The case where $\omega_1(\bold c)$ is indefinite.
\endhead
    Each indefinite quadratic form in $\Bbb R^2$ has a basis of two 
linearly independent homogeneous vectors $\bold c_1$ and $\bold c_2$, 
i\.\,e\. two vectors in $\Bbb R^2$ such that 
$$
\xalignat 2
&\hskip -2em
\omega_1(\bold c_1)=0,
&&\omega_1(\bold c_2)=0.
\mytag{4.1} 
\endxalignat
$$
In physics the homogeneous vectors $\bold c_1$ and $\bold c_2$ in \mythetag{4.1} 
are called light vectors or light cone vectors (see \mycite{77}). Choosing
$$
\xalignat 2
&\bold c_1=\Vmatrix c_{11}\\ 
\vspace{1ex}c_{21}\!\endVmatrix,
&&\bold c_2=\Vmatrix c_{12}\\
\vspace{1ex}c_{22}\!\endVmatrix,
\endxalignat
$$
we construct a non-degenerate matrix $C$ in \mythetag{2.17} and apply it in
\mythetag{2.13}. For such a matrix both equalities \mythetag{2.14} and
\mythetag{2.15} are fulfilled. Hence two of the three equalities \mythetag{2.12}
are fulfilled. The third equality \mythetag{2.12} is not fulfilled since 
otherwise we would return to the previous case where $\omega_1(\bold c)$ is
zero: 
$$
\xalignat 3
&\det\Vmatrix a_{20} & a_{11}\!\\
\vspace{1ex}
b_{20} & b_{11}\!\endVmatrix=0,
&&\det\Vmatrix a_{11} & a_{02}\!\\
\vspace{1ex}
b_{11} & b_{02}\!\endVmatrix=0,
&&\det\Vmatrix a_{20} & a_{02}\!\\
\vspace{1ex}
b_{20} & b_{02}\!\endVmatrix\neq 0.
\qquad\quad
\mytag{4.2}
\endxalignat
$$
The relationships mean that the first two of the three vectors 
$$
\xalignat 3
&\hskip -2em
\Vmatrix a_{20}\!\\ \vspace{1ex} b_{20}\!\endVmatrix,
&&\Vmatrix a_{02}\!\\ \vspace{1ex} b_{02}\!\endVmatrix,
&&\Vmatrix a_{11}\!\\ \vspace{1ex} b_{11}\!\endVmatrix
\mytag{4.3}
\endxalignat
$$
in \mythetag{4.3} are linearly independent, while the third vector 
\mythetag{4.3} belongs to the span of each of the first two. This fact
implies $a_{11}=0$ and $b_{11}=0$. It means that by applying some properly 
chosen linear transformation \mythetag{2.13}, we can bring our 
quadratic transformation \mythetag{2.1} to the following form: 
$$
\hskip -3em
\aligned
&\tilde p=a_{20}\,p^2+a_{02}\,q^2+2\,a_{10}\,p
+2\,a_{01}\,q+a_{00},\\
&\tilde q=b_{20}\,p^2+b_{02}\,q^2+2\,b_{10}\,p
+2\,b_{01}\,q+b_{00}.
\endaligned
\mytag{4.4}
$$\par
     Now lets consider another linear transformation similar to 
\mythetag{2.13}:
$$
\xalignat 2
&\hskip -2em
\hat p=d_{11}\,\tilde p+d_{12}\,\tilde q,
&&\hat q=d_{21}\,\tilde p+d_{22}\,\tilde q.
\mytag{4.5}
\endxalignat
$$
The transformations \mythetag{2.3}, \mythetag{2.4}, and \mythetag{2.10}
are special instances of the transformation \mythetag{4.5}. Since the third
determinant in \mythetag{4.2} is nonzero, choosing properly the coefficients
of \mythetag{4.5} and then applying \mythetag{4.5} to \mythetag{4.4}, we 
can bring \mythetag{4.4} to 
$$
\hskip -3em
\aligned
&\tilde p=p^2+2\,a_{10}\,p+2\,a_{01}\,q+a_{00},\\
&\tilde q=q^2+2\,b_{10}\,p+2\,b_{01}\,q+b_{00}.
\endaligned
\mytag{4.6}
$$
Using origin shifts we can further simplify \mythetag{4.6} bringing it to 
$$
\hskip -3em
\aligned
&\tilde p=p^2+2\,a_{01}\,q,\\
&\tilde q=q^2+2\,b_{10}\,p.
\endaligned
\mytag{4.7}
$$\par
     The coefficient $2\,b_{10}$ in \mythetag{4.7} is nonzero. Indeed, 
otherwise we would have $\tilde q=q^2$. The equality $\tilde q=q^2$ cannot
be resolved with respect to $q$ if $\tilde q<0$ thus contradicting the 
invertibility of the transformation \mythetag{4.7}. Similar arguments apply
to $2\,a_{01}$ in \mythetag{4.7}. Therefore the following conditions are 
fulfilled:
$$
\xalignat 2
&\hskip -2em
a_{01}\neq 0,
&&b_{10}\neq 0.
\mytag{4.8}
\endxalignat
$$
Relying on \mythetag{4.8}, we consider the scaling transformations
$$
\xalignat 2
&\hskip -2em
p=\alpha\,\check p,
&&q=\beta\,\check q,
\mytag{4.9}\\
&\hskip -2em
\hat p=\alpha^2\,\tilde p,
&&\hat q=\beta^2\,\tilde q,
\mytag{4.10}\\
\endxalignat
$$
where $\alpha$ and $\beta$ are given by the formulas
$$
\xalignat 2
&\hskip -2em
\alpha=2\,\root 3\of{a_{01}^{\kern 2pt 2}\,b_{10}},
&&\beta=2\,\root 3\of{a_{01}\,b_{10}^{\kern 2pt 2}}.
\mytag{4.11}
\endxalignat
$$
The transformations \mythetag{4.9} and \mythetag{4.10} are special
instances of \mythetag{2.13} and \mythetag{4.5}. Applying them 
to \mythetag{4.7} and taking into account \mythetag{4.11}, we derive
$$
\xalignat 2
&\hskip -2em
\tilde p=p^2+q,
&&\tilde q=q^2+p.
\mytag{4.12}
\endxalignat
$$\par
     The second equality \mythetag{4.12} can be resolved with respect to 
the variable $p$. It yields ${p=\tilde q-q^2}$. Substituting $p=\tilde q-q^2$
into the first equation \mythetag{4.12}, \pagebreak we derive the following 
quartic equation with respect to the variable $q$:
$$
\hskip -2em
q^4-2\,\tilde q\,q^2+q+\tilde q^2-\tilde p=0.
\mytag{4.13}
$$ 
The invertibility of the transformation \mythetag{4.12} means that for any 
real values of $\tilde p$ and $\tilde q$ the quartic equation \mythetag{4.13}
should have exactly one real root. Fortunately there is a criterion for a
quartic equation with real coefficients to have exactly one real root. This
criterion is given below by Theorem~\mythetheorem{A.1} in Appendix A.\par
     Theorem~\mythetheorem{A.1} says that the discriminant of the equation 
\mythetag{4.13} should be zero:
$$
\hskip -2em
D_4=0.
\mytag{4.14}
$$
The rest is to calculate the discriminant $D_4$ in \mythetag{4.14} explicitly. 
Applying the formula \mythetag{A.7} to the coefficients of the equation
\mythetag{4.13}, we derive: 
$$
\hskip -2em
D_4=-256\,\tilde p^3+256\,\tilde q^2\,\tilde p^2+288\,\tilde q\,\tilde p
-256\,\tilde q^3-27. 
\mytag{4.15}
$$
The discriminant \mythetag{4.15} depends on $\tilde p$ and $\tilde q$. It is not 
identically zero. Therefore the equality \mythetag{4.14} cannot be fulfilled for
all real values of $\tilde p$ and $\tilde q$. This result is formulated
as the following theorem. 
\mytheorem{4.1} There is no invertible quadratic transformation in $\Bbb R^2$
whose associated quadratic form $\omega_1(\bold c)$ in \mythetag{2.18} is
indefinite. 
\endproclaim
\head
5. The case where $\omega_1(\bold c)$ is semi-definite
\endhead
     Each semi-definite quadratic form in $\Bbb R^2$ up to a scalar factor 
has only one homogeneous vector $\bold c_1$, i\.\,e\. a vector in $\Bbb R^2$ 
such that $\omega_1(\bold c_1)=0$. Choosing $\bold c_1$ for one of the two 
columns of the matrix \mythetag{2.17}, we can satisfy only one of the 
two relationships \mythetag{2.14} or \mythetag{2.15}. Let's choose $\bold c_1$
for the first column of the matrix \mythetag{2.17}. Then the relationship 
\mythetag{2.14} is fulfilled, while the relationship \mythetag{2.15} is not
fulfilled. Upon applying the transformation \mythetag{2.13} to \mythetag{2.1}
we shall find that the first equality \mythetag{2.12} is fulfilled, while
the second equality \mythetag{2.12} is broken: 
$$
\xalignat 2
&\det\Vmatrix a_{20} & a_{11}\!\\
\vspace{1ex}
b_{20} & b_{11}\!\endVmatrix=0,
&&\det\Vmatrix a_{11} & a_{02}\!\\
\vspace{1ex}
b_{11} & b_{02}\!\endVmatrix\neq 0.
\qquad\quad
\mytag{5.1}
\endxalignat
$$
The second relationship \mythetag{5.1} means that the vectors
$$
\xalignat 2
&\Vmatrix a_{02}\!\\ \vspace{1ex} b_{02}\!\endVmatrix,
&&\Vmatrix a_{11}\!\\ \vspace{1ex} b_{11}\!\endVmatrix
\endxalignat
$$
are linearly independent. Hence they are nonzero. This means that
$a_{11}\neq 0$ or $b_{11}\neq 0$. Applying the transformation 
\mythetag{2.4} if needed, we can assume that $a_{11}\neq 0$. Hence 
we can apply \mythetag{2.3} and derive $b_{11}=0$ after that. The
first relationship \mythetag{5.1} means that we shall derive 
$b_{20}=0$ along with $b_{11}=0$, i\.\,e\. the relationships 
\mythetag{2.6} will be fulfilled. Once they are fulfilled, we find 
that \mythetag{2.1} is transformed to 
$$
\hskip -3em
\aligned
&\tilde p=a_{20}\,p^2+2\,a_{11}\,p\,q+a_{02}\,q^2+2\,a_{10}\,p
+2\,a_{01}\,q+a_{00},\\
&\tilde q=b_{02}\,q^2+2\,b_{10}\,p
+2\,b_{01}\,q+b_{00}.
\endaligned
\mytag{5.2}
$$
The coefficient $b_{02}$ in \mythetag{5.2} is nonzero, since otherwise
we return to the previous case where the quadratic form $\omega_1(\bold c)$
is zero: 
$$
\hskip -2em
b_{02}\neq 0.
\mytag{5.3}
$$\par     
Using \mythetag{5.3} and composing \mythetag{5.2} with the linear 
transformation
$$
\xalignat 2
&\hat p=p-\frac{a_{02}}{b_{02}}\,q,
&&\hat q=q, 
\endxalignat
$$
we can bring the transformation \mythetag{5.2} to the form with $a_{02}=0$:
$$
\hskip -3em
\aligned
&\tilde p=a_{20}\,p^2+2\,a_{11}\,p\,q+2\,a_{10}\,p
+2\,a_{01}\,q+a_{00},\\
&\tilde q=b_{02}\,q^2+2\,b_{10}\,p
+2\,b_{01}\,q+b_{00}.
\endaligned
\mytag{5.4}
$$\par
     Let's apply \mythetag{5.4} to \mythetag{2.18}. The first determinant in 
\mythetag{2.18} vanishes. The second and the third determinants for
\mythetag{5.4} are calculated explicitly:
$$
\xalignat 3
&\vmatrix a_{20} & a_{11}\!\\
\vspace{1ex}
b_{20} & b_{11}\!\endvmatrix=0,
&&\vmatrix a_{20} & a_{02}\!\\
\vspace{1ex}
b_{20} & b_{02}\!\endvmatrix=a_{20}\,b_{02},
&&\vmatrix a_{11} & a_{02}\!\\
\vspace{1ex}
b_{11} & b_{02}\!\endvmatrix=a_{11}\,b_{02}.
\qquad
\mytag{5.5}
\endxalignat
$$
Due to \mythetag{5.5} and \mythetag{2.18} the quadratic form 
$\omega_1(\bold c)$ is presented by the matrix 
$$
\xalignat 2
&\hskip -2em
\Omega_1=\Vmatrix 0 & \dfrac{a_{20}\,b_{02}}{2}\\
\vspace{2ex}
\dfrac{a_{20}\,b_{02}}{2} & a_{11}\,b_{02}\endVmatrix,
&&\det \Omega_1=-\frac{a_{20}^{\kern 2pt 2}\,b_{02}^{\kern 2pt 2}}{4}.
\mytag{5.6}
\endxalignat
$$
It is known that a quadratic form in $\Bbb R^2$ is semi-definite if and
only if its matrix is nonzero, but the determinant of its matrix is zero
(see \mycite{75}). Applying this fact to \mythetag{5.6} and taking into
account \mythetag{5.3}, we derive 
$$
\xalignat 2
&\hskip -2em
a_{20}=0,
&&a_{11}\neq 0.
\mytag{5.7}
\endxalignat
$$
Due to \mythetag{5.7} the transformation \mythetag{5.4} looks like
$$
\hskip -3em
\aligned
&\tilde p=2\,a_{11}\,p\,q+2\,a_{10}\,p
+2\,a_{01}\,q+a_{00},\\
&\tilde q=b_{02}\,q^2+2\,b_{10}\,p
+2\,b_{01}\,q+b_{00}.
\endaligned
\mytag{5.8}
$$
By means of scaling transformations and origin shifts we can bring 
\mythetag{5.8} to 
$$
\hskip -3em
\aligned
&\tilde p=p\,q+2\,a_{10}\,p
+2\,a_{01}\,q+a_{00},\\
&\tilde q=q^2+2\,b_{10}\,p.
\endaligned
\mytag{5.9}
$$\par
     Using the same arguments as in deriving \mythetag{5.9}, we 
find that $\,b_{10}\neq 0$ in \mythetag{5.9}. Hence, applying proper
scaling transformation, we can further simplify \mythetag{5.9}:
$$
\hskip -3em
\aligned
&\tilde p=p\,q+2\,a_{10}\,p
+2\,a_{01}\,q+a_{00},\\
&\tilde q=q^2+p.
\endaligned
\mytag{5.10}
$$
Then by means of origin shifts in $p$, $\tilde q$, and $\tilde p$ 
we can annul $a_{01}$ and $a_{00}$ in \mythetag{5.10}:
$$
\hskip -3em
\aligned
&\tilde p=p\,q+2\,a_{10}\,p,\\
&\tilde q=q^2+p.
\endaligned
\mytag{5.11}
$$
The coefficient $2\,a_{10}$ in \mythetag{5.11} can be either zero or nonzero.
In both cases we can resolve the second equality \mythetag{5.11} with respect
to $p$:
$$
\hskip -2em
p=\tilde q-q^2.
\mytag{5.12}
$$
Substituting \mythetag{5.12} back to the first equality \mythetag{5.11}, we
derive the cubic equation 
$$
\hskip -2em
q^3+2\,a_{10}\,q^2-\tilde q\,q-2\,a_{10}\,\tilde q+\tilde p=0
\mytag{5.13}
$$
for the variable $q$. The invertibility of the transformation 
\mythetag{5.11} means that for any real values of $\tilde p$ and 
$\tilde q$ the cubic equation \mythetag{4.14} should have exactly one 
real root. Fortunately there is a criterion for a general cubic equation 
with real coefficients 
$$
\hskip -2em
q^3+\alpha_1\,q^2+\alpha_2\,q+\alpha_3=0
\mytag{5.14}
$$
to have exactly one real root. It is formulated through its discriminant
$$
\hskip -2em
D_3=-27\,\alpha_3^2+18\,\alpha_3\,\alpha_1\,\alpha_2+\alpha_1^2\,\alpha_2^2
-4\,\alpha_1^3\,\alpha_3-4\,\alpha_2^3.
\mytag{5.15}
$$
\mytheorem{5.1} A cubic equation with  real coefficients \mythetag{5.14}
has exactly one simple real root if and only if its discriminant \mythetag{5.15}
is negative: $D_3<0$.
\endproclaim
Theorem~\mythetheorem{5.1} is immediate from the following theorem.
\mytheorem{5.2} A cubic equation with  real coefficients \mythetag{5.14}
has three distinct real roots if and only if its discriminant \mythetag{5.15}
is positive: $D_3>0$.
\endproclaim
     Indeed, the case $D_3=0$ is trivial. In this case the equation 
\mythetag{5.14} has one real root (of multiplicity $3$) or two real roots
(one of which is simple and the other is double). As for 
Theorem~\mythetheorem{5.2}, it is well-known. Its proof can be found in 
\mycite{78} (see Theorem~3.1 over there).\par
     Now, applying Theorem~\mythetheorem{5.1} to the equation 
\mythetag{5.13} and taking into account the case $D_3=0$ with one triple
root, we write the inequality
$$
\hskip -2em
D_3\leqslant 0. 
\mytag{5.16}
$$
The rest is to calculate the discriminant $D_3$ of the equation
\mythetag{5.13} explicitly:
$$
\hskip -2em
D_3=4\,\tilde q^3-32\,a_{10}^{\kern 2pt 2}\,\tilde q^2+8\,a_{10}
\,(9\,\tilde p+8\,a_{10}^{\kern 2pt 3})\,\tilde q
-\tilde p\,(32\,a_{\kern 2pt 10}^3+27\,\tilde p).
\mytag{5.17}
$$
It is easy to see that the discriminant \mythetag{5.17} is a cubic 
polynomial with respect to the variable $\tilde q$ with the non-vanishing
leading coefficient $4$. Such a polynomial takes values of both signs ---
positive and negative. Therefore the inequality \mythetag{5.16} cannot be
fulfilled for all $\tilde p$ and $\tilde q$. 
\mytheorem{5.3} There is no invertible quadratic transformation in $\Bbb R^2$
whose associated quadratic form $\omega_1(\bold c)$ in \mythetag{2.18} is
semi-definite. 
\endproclaim
\head
5. The case where $\omega_1(\bold c)$ is definite.
\endhead
\adjustfootnotemark{-1}
     If $\omega_1(\bold c)$ in \mythetag{2.18} is definite, none of the 
equalities \mythetag{2.14} and \mythetag{2.15} can be fulfilled. Therefore
we proceed to the equality \mythetag{2.16}. Looking at it, we define the second 
quadratic form $\omega_2(\bold c)$ associated with a quadratic transformation 
of the form \mythetag{2.1}. Unlike $\omega_1(\bold c)$, the form $\omega_2(\bold c)$
is a quadratic form\footnotemark\ in the four-dimensional vector space $\Bbb R^4$. It is given by
the following formula: 
$$
\aligned
&\omega_2(\bold c)=c_1\,\vmatrix a_{20} & a_{11}\!\\
\vspace{1ex}
b_{20} & b_{11}\!\endvmatrix\,c_3+c_1\,\vmatrix a_{20} & a_{02}\!\\
\vspace{1ex}
b_{20} & b_{02}\!\endvmatrix\,c_4\,+\\
&\kern 7.7em+\,c_2\,\vmatrix a_{20} & a_{02}\!\\
\vspace{1ex}
b_{20} & b_{02}\!\endvmatrix\,c_3+c_2\,\vmatrix a_{11} & a_{02}\!\\
\vspace{1ex}
b_{11} & b_{02}\!\endvmatrix\,c_4.
\endaligned
\mytag{6.1}
$$
Here $c_1$, $c_2$, $c_3$, $c_4$ are the coordinates of a vector 
$\bold c\in\Bbb R^4$.\par
\footnotetext{\ Regular quadratic forms are tensors of the type $(0,2)$. Like
$\omega_1(\bold c)$, the quadratic form $\omega_2(\bold c)$ is a pseudo-tensor 
of the type $(0,2)$ in the sense of Definition~2.1 in \mycite{76}.} 
     The determinant of the matrix $C$ from \mythetag{2.17} produces another
quadratic form in $R^4$. It is given by the following formula:
$$
\hskip -2em
\omega_3(\bold c)=2\,c_1\,c_4-\,2\,c_2\,c_3=2\,\det C. 
\mytag{6.2}
$$
Let's denote through $2\,\alpha$, $2\,\beta$, and $2\,\gamma$ the determinants in 
\mythetag{6.1} and \mythetag{2.18}:
$$
\xalignat 3
&\hskip -2em
2\,\alpha=\vmatrix a_{20} & a_{11}\!\\
\vspace{1ex}
b_{20} & b_{11}\!\endvmatrix,
&&2\,\beta=
\vmatrix a_{20} & a_{02}\!\\
\vspace{1ex}
b_{20} & b_{02}\!\endvmatrix,
&&2\,\gamma=\vmatrix a_{11} & a_{02}\!\\
\vspace{1ex}
b_{11} & b_{02}\!\endvmatrix.
\qquad
\mytag{6.3}
\endxalignat
$$
Then, like in \mythetag{5.6}, we can write the matrix presentation of 
\mythetag{2.18}:
$$
\Omega_1=\Vmatrix 2\,\alpha & \beta\\
\vspace{2ex}
\beta & 2\,\gamma\endVmatrix.
\mytag{6.4}
$$
The forms \mythetag{6.1} and \mythetag{6.2} have the following matrix
presentations
$$
\xalignat 2
&\Omega_2=\Vmatrix 0 & 0 & \alpha & \beta\\
\vspace{0.5ex}
0 & 0 & \beta & \gamma\\
\vspace{0.5ex}
\alpha & \beta & 0 & 0\\
\vspace{0.5ex}
\beta & \gamma & 0 & 0
\endVmatrix,
&&\Omega_3=\Vmatrix 0 & 0 & 0 & 1\\
\vspace{0.5ex}
0 & 0 & \kern -3pt -1 & 0\\
\vspace{0.5ex}
0 & \kern -3pt -1 & 0 & 0\\
\vspace{0.5ex}
1 & 0 & 0 & 0
\endVmatrix.
\endxalignat
$$
It is known that a quadratic form in $\Bbb R^2$ is definite if and only if
the diagonal elements of its matrix are nonzero and of the same sign, and
if the determinant of its matrix is positive (see Silvester's criterion in
\mycite{75}). For \mythetag{6.4} this yields 
$$
\xalignat 2
&\hskip -2em
\alpha\,\gamma>0,
&&4\,\alpha\,\gamma-\beta^2>0.
\mytag{6.5}
\endxalignat
$$\par
     The form $\omega_2(\bold c)$ in \mythetag{6.1} cannot be zero since otherwise
the form $\omega_1(\bold c)$ in \mythetag{2.18} would be zero thus returning us
to one of the previous cases. Note that there are no squares of $c_1$, $c_2$, $c_3$,
$c_4$ in \mythetag{6.1}. Hence the form $\omega_2(\bold c)$ is either indefinite
or semi-definite. In both cases there is at least one vector $\bold c_1\in\Bbb R^4$ 
such that 
$$ 
\hskip -2em
\omega_2(\bold c_1)=0.
\mytag{6.6}
$$
Apart from \mythetag{6.6}, we need to fulfill the other condition
$$
\hskip -2em
\omega_3(\bold c_1)\neq 0.
\mytag{6.7}
$$
Therefore we shall calculate the components of $\bold c_1$ explicitly. Using the
notations \mythetag{6.3} and applying them to \mythetag{6.1}, we write 
\mythetag{6.6} as 
$$
\gather
\hskip -2em
\alpha\,c_1\,c_3+\beta\,c_1\,c_4+\beta\,c_2\,c_3+\gamma\,c_2\,c_4=0,\\
\vspace{-2.5ex}
\intertext{or equivalently}
\vspace{-2.5ex}
\hskip -2em
c_1\,(\alpha\,c_3+\beta\,c_4)+c_2(\beta\,c_3+\gamma\,c_4)=0.
\mytag{6.8}
\endgather
$$
Let's choose the following values for $c_1$, $c_2$, $c_3$, $c_4$:
$$
\xalignat 4
&\hskip -2em
c_1=\gamma,
&&c_2=-3\,\beta,
&&c_3=-2\,\beta,
&&c_4=\frac{3\,\beta^2-\alpha\,\gamma}{\gamma}. 
\qquad
\mytag{6.9}
\endxalignat
$$
It is easy to see that \mythetag{6.9} is a solution of the equation 
\mythetag{6.8}. Substituting \mythetag{6.9} into \mythetag{6.2} and
taking into account \mythetag{6.5}, we derive 
$$
\hskip -2em
\omega_3(\bold c_1)=-2\,(\alpha\,\gamma+3\,\beta^2)
<-\frac{13}{2}\,\beta^2<0.
\mytag{6.10}
$$
The inequality \mythetag{6.10} means that the inequality \mythetag{6.7} is
fulfilled.\par
    Now we use the quantities \mythetag{6.9} as the components of the matrix 
$C$ in \mythetag{2.17}:
$$
\xalignat 2
&\hskip -2em
c_{11}=c_1=\gamma,
&&c_{21}=c_2=-3\,\beta,\\
\vspace{-1.5ex}
\mytag{6.11}\\
\vspace{-1.5ex}
&\hskip -2em
c_{12}=c_3=-2\,\beta,
&&c_{22}=c_4=\frac{3\,\beta^2-\alpha\,\gamma}{\gamma}. 
\endxalignat
$$
The inequality \mythetag{6.10} implying \mythetag{6.7} means that 
$\det C\neq 0$. Therefore we can use the matrix $C$ with the components 
\mythetag{6.11} in \mythetag{2.13}. Upon applying the transformation 
\mythetag{2.13} to \mythetag{2.1} we shall find that the third equality 
is \mythetag{2.12} fulfilled. As for the first two equalities \mythetag{2.12},
they are broken since otherwise the form $\omega_1(\bold c)$ would not be
definite thus returning us to one of the previous cases:
$$
\xalignat 3
&\det\Vmatrix a_{20} & a_{11}\!\\
\vspace{1ex}
b_{20} & b_{11}\!\endVmatrix\neq 0,
&&\det\Vmatrix a_{11} & a_{02}\!\\
\vspace{1ex}
b_{11} & b_{02}\!\endVmatrix\neq 0,
&&\det\Vmatrix a_{20} & a_{02}\!\\
\vspace{1ex}
b_{20} & b_{02}\!\endVmatrix=0.
\qquad\quad
\mytag{6.12}
\endxalignat
$$\par
     The first relationship \mythetag{6.12} means that the vectors
$$
\xalignat 2
&\hskip -2em
\Vmatrix a_{20}\!\\ \vspace{1ex} b_{20}\!\endVmatrix,
&&\Vmatrix a_{11}\!\\ \vspace{1ex} b_{11}\!\endVmatrix
\mytag{6.13}
\endxalignat
$$
are linearly independent. Hence they are nonzero. This means that
$a_{20}\neq 0$ or $b_{20}\neq 0$. Applying the transformation 
\mythetag{2.4} if needed, we can assume that $a_{20}\neq 0$. Hence 
we can apply \mythetag{2.10} and derive $b_{20}=0$ after that. The
last relationship \mythetag{6.13} means that we shall derive 
$b_{02}=0$ along with $b_{21}=0$, i\.\,e\. the relationships 
\mythetag{2.11} will be fulfilled. Once they are fulfilled, we find 
that \mythetag{2.1} is transformed to 
$$
\hskip -3em
\aligned
&\tilde p=a_{20}\,p^2+2\,a_{11}\,p\,q+a_{02}\,q^2+2\,a_{10}\,p
+2\,a_{01}\,q+a_{00},\\
&\tilde q=2\,b_{11}\,p\,q+2\,b_{10}\,p+2\,b_{01}\,q+b_{00}.
\endaligned
\mytag{6.14}
$$
The coefficient $b_{11}$ in \mythetag{6.14} is nonzero, since otherwise
we return to the previous case where the quadratic form $\omega_1(\bold c)$
is zero: 
$$
\hskip -2em
b_{11}\neq 0.
\mytag{6.15}
$$\par     
     Using \mythetag{6.15} and composing \mythetag{6.14} with a properly
chosen linear transformation, we can bring the transformation \mythetag{6.14} 
to the form with $a_{11}=0$:
$$
\hskip -3em
\aligned
&\tilde p=a_{20}\,p^2+a_{02}\,q^2+2\,a_{10}\,p+2\,a_{01}\,q+a_{00},\\
&\tilde q=2\,b_{11}\,p\,q+2\,b_{10}\,p+2\,b_{01}\,q+b_{00}.
\endaligned
\mytag{6.16}
$$ 
Let's apply \mythetag{6.16} to \mythetag{2.18}. The second determinant in 
\mythetag{2.18} vanishes. The first and the third determinants for
\mythetag{6.16} are calculated explicitly:
$$
\xalignat 3
&\vmatrix a_{20} & a_{11}\!\\
\vspace{1ex}
b_{20} & b_{11}\!\endvmatrix=a_{20}\,b_{11},
&&\vmatrix a_{20} & a_{02}\!\\
\vspace{1ex}
b_{20} & b_{02}\!\endvmatrix=0,
&&\vmatrix a_{11} & a_{02}\!\\
\vspace{1ex}
b_{11} & b_{02}\!\endvmatrix=-a_{02}\,b_{11}.
\qquad
\mytag{6.17}
\endxalignat
$$
Due to \mythetag{6.17} and \mythetag{2.18} the quadratic form 
$\omega_1(\bold c)$ is presented by the matrix 
$$
\xalignat 2
&\hskip -2em
\Omega_1=\Vmatrix a_{20}\,b_{11} & 0\\
\vspace{2ex}
0 & -a_{02}\,b_{11}\endVmatrix,
&&\det \Omega_1=-a_{20}\,a_{02}\,b_{11}^{\kern 2pt 2}.
\mytag{6.18}
\endxalignat
$$
It is known that a quadratic form in $\Bbb R^2$ is definite if and
only if the diagonal elements of its matrix are nonzero and of the same sign, 
and if the determinant of its matrix is positive (see Silvester's criterion in
\mycite{75}). Applying this fact to \mythetag{6.18} and taking into
account \mythetag{6.15}, we derive 
$$
\xalignat 3
&\hskip -2em
a_{20}\neq 0,
&&a_{02}\neq 0,
&&a_{20}\,a_{02}<0.
\mytag{6.19}
\endxalignat
$$\par 
     Using the inequalities \mythetag{6.19}, by means of scaling transformations 
and origin shifts we can bring \mythetag{6.16} to the following form:
$$
\xalignat 2
&\hskip -2em
\tilde p=p^2-q^2+2\,a_{10}\,p+2\,a_{01}\,q,
&&\tilde q=p\,q.
\mytag{6.20}
\endxalignat
$$ 
The second equality \mythetag{6.20} is linear with respect to $p$. We can resolve 
it as 
$$
\hskip -2em
p=\frac{\tilde q}{q}\text{\ \ if \ }q\neq 0. 
\mytag{6.21}
$$
Let's substitute \mythetag{6.21} into the first equality \mythetag{6.20}. Removing
denominators, we derive the following quartic equation with respect to $q$:
$$
\hskip -2em
q^4+2\,a_{01}\,q^3-\tilde p\,q^2+2\,a_{10}\,\tilde q\,q+\tilde q^2=0.
\mytag{6.22}
$$
If $\tilde q\neq 0$, the exceptional value $q=0$ is not a root of the quartic
equation \mythetag{6.22}. For all $\tilde p$ and $\tilde q\neq 0$ the invertibility 
of the transformation \mythetag{6.5} means that the quartic equation \mythetag{6.22} 
has exactly one real root. Fortunately there is a criterion for a quartic equation 
with real coefficients to have exactly one real root. This criterion is given below 
by Theorem~\mythetheorem{A.1} in Appendix A.\par
     Theorem~\mythetheorem{A.1} says that the discriminant of the equation 
\mythetag{6.22} should be zero:
$$
\hskip -2em
D_4=0.
\mytag{6.23}
$$
The rest is to calculate the discriminant $D_4$ in \mythetag{6.23} explicitly. 
Applying the formula \mythetag{A.7} to the coefficients of the equation
\mythetag{6.22}, we derive: 
$$
\hskip -2em
\gathered
D_4=256\,\tilde q^6
-768\,a_{01}\,a_{10}\,\tilde q^5
-(576\,\tilde p\,a_{01}^2+576\,a_{10}^2\,\tilde p\,+\\
+\,432\,a_{01}^4+96\,a_{01}^2\,a_{10}^2+128\,\tilde p^2
+432\,a_{10}^4)\,\tilde q^4-(288\,a_{01}^3\,a_{10}\,\tilde p\,+\\
+\,320\,a_{01}\,a_{10}\,\tilde p^2+256\,a_{01}^3\,a_{10}^3
+288\,a_{01}\,a_{10}^3\,\tilde p)\,\tilde q^3+(16\,\tilde p^4\,+\\
+\,16\,\tilde p^3\,a_{10}^2+16\,\tilde p^2\,a_{01}^2\,a_{10}^2
+16\,\tilde p^3\,a_{01}^2)\,\tilde q^2.
\endgathered
\mytag{6.24}
$$
As we see, the discriminant \mythetag{6.24} is a polynomial of sixth
degree. For $\tilde q\neq 0$ this polynomial is not identically zero. 
Therefore the equality \mythetag{6.23} cannot be fulfilled for all
$\tilde p$ and $\tilde q\neq 0$. This result leads to the following
theorem. 
\mytheorem{6.1} There is no invertible quadratic transformation in $\Bbb R^2$
whose associated quadratic form $\omega_1(\bold c)$ in \mythetag{2.18} is
definite. 
\endproclaim
\head
6. Conclusions.
\endhead
\mydefinition{6.1} Two quadratic mappings $f_1\!:\,\Bbb R^2\to \Bbb R^2$ and
$f_2\!:\,\Bbb R^2\to \Bbb R^2$ are called equivalent if there are two invertible 
linear mappings $\varphi_1\!:\,\Bbb R^2\to \Bbb R^2$ and
$\varphi_2\!:\,\Bbb R^2\to \Bbb R^2$ such that $\varphi_1\compos f_1=
f_2\compos\varphi_2$.
\enddefinition
    The main result of the present paper consist in subdividing all potentially
invertible quadratic transformations of the real plane $\Bbb R^2$ into four
groups and in finding up to the equivalence introduced in 
Definition~\mythedefinition{6.1} some pre-canonical presentations for quadratic 
transformations within these groups. These presentations are given by the 
formulas \mythetag{3.5}, \mythetag{4.12}, \mythetag{5.11}, and \mythetag{6.20}. 
Theorems~\mythetheorem{3.1}, \mythetheorem{4.1}, \mythetheorem{5.3}, and
\mythetheorem{6.1} show that only transformations of the first of the four
groups are actually invertible.\par 
     Some examples of quadratic transformations of $\Bbb R^2$ were studied
in \mycite{79} as discrete dynamical systems. Some prospects of applying 
quadratic transformations to the perfect cuboid problem are discussed in the 
introductory section of this paper.\par
\newpage
\head
Appendix A.\\
Quartic~polynomials~with~exactly~one~real~root.
\endhead
\rightheadtext{Appendix A. Quartic polynomials \dots}
\leftheadtext{Ruslan Sharipov}
\parshape 20 0pt 360pt 0pt 360pt 0pt 360pt 0pt 360pt 0pt 360pt 0pt 360pt
170pt 190pt 170pt 190pt 170pt 190pt 170pt 190pt 170pt 190pt 170pt 190pt
170pt 190pt 170pt 190pt 170pt 190pt 170pt 190pt 170pt 190pt 170pt 190pt
170pt 190pt 0pt 360pt 
    Let $P_4(x)$ be a real quartic polynomial. Without loss of generality we can 
assume that $P_4(x)$ is monic, i\.\,e\. its leading coefficient is equal to one:
$$
\hskip -2em
P_4(x)=x^4+a_1\,x^3+a_2\,x^2+a_3\,x+a_4.
\mytag{A.1}
$$
Our goal is to find a necessary and sufficient condition for the polynomial 
\mythetag{A.1} to have a unique real root. 
\vadjust{\vskip 130pt\hbox to 0pt{\kern -10pt
\includegraphics{quadro_trans_a1.eps}\hss}\vskip -130pt} The graph of such 
a polynomial is shown in Fig\.~A.1. Since $P_4(x)\to +\infty$ as $x\to\pm\infty$ the
unique root $x_0$ of the polynomial \mythetag{A.1} should be a root of multiplicity
at least two. Hence we have 
$$
\hskip -2em
P_4(x)=(x-x_0)^2\,P_2(x). 
\mytag{A.2}
$$
Here $P_2(x)$ is a quadratic polynomial with real coefficients. We write it as
$$
\hskip -2em
P_2(x)=x^2+b_1\,x+b_2.
\mytag{A.3}
$$
Substituting \mythetag{A.3} into \mythetag{A.2} and expanding the resulting 
expression we can express the coefficients of the polynomial $P_2(x)$ through
$a_1$, $a_2$, and $x_0$:
$$
\xalignat 2
&\hskip -2em
b_1=a_1+2\,x_0,
&&b_2=a_2-x_0^2+2\,b_1\,x_0.
\mytag{A.4}
\endxalignat
$$ 
Moreover, we obtain two formulas expressing $a_3$ and $a_4$
through $b_1$, $b_2$, and $x_0$:
$$
\xalignat 2
&\hskip -2em
a_4=b_2\,x_0^2,
&&a_3=b_1\,x_0^2-2\,b_2\,x_0.
\mytag{A.5}
\endxalignat
$$\par
     The root $x_0$ in \mythetag{A.2} is not simple. This means that the
discriminant of the quartic polynomial $P_4(x)$ in \mythetag{A.1} is 
equal to zero:
$$
\hskip -2em
D_4=0.
\mytag{A.6}
$$
The explicit formula for the discriminant $D_4$ in \mythetag{A.6} looks like
$$
\gathered
D_4=18\,a_1^3\,a_3\,a_2\,a_4+256\,a_4^3-6\,a_1^2\,a_3^2\,a_4-192\,a_1\,a_3\,a_4^2
+18\,a_1\,a_3^3\,a_2\,+\\
+\,144\,a_2\,a_1^2\,a_4^2+a_2^2\,a_1^2\,a_3^2-4\,a_2^3\,a_1^2\,a_4
+144\,a_4\,a_3^2\,a_2-4\,a_1^3\,a_3^3-27\,a_3^4\,-\\
-\,128\,a_2^2\,a_4^2+16\,a_2^4\,a_4-4\,a_2^3\,a_3^2-27\,a_1^4\,a_4^2
-80\,a_1\,a_3\,a_2^2\,a_4.
\endgathered
\quad
\mytag{A.7}
$$
The equality \mythetag{A.6} is a necessary condition for the real quartic 
polynomial \mythetag{A.1} to have a unique real root. But this condition 
is not sufficient.\par 
     The quadratic polynomial $P_2(x)$ in \mythetag{A.2} should have at most 
one real root coinciding with $x_0$. Therefore its discriminant $D_2$ should
be non-positive: 
$$
\hskip -2em
D_2\leqslant 0. 
\mytag{A.8}
$$
The discriminant of the polynomial $P_2(x)$ is given by the formula
$$
\hskip -2em
D_2=b_1^2-4\,b_2.
\mytag{A.9}
$$
Substituting \mythetag{A.4} into \mythetag{A.9}, we derive the following formula:
$$
\hskip -2em
D_2=a_1^2-4\,a_2-8\,x_0^2-4\,x_0\,a_1. 
\mytag{A.10}
$$\par
     The formula \mythetag{A.10} comprises the root $x_0$. In order to write the
inequality \mythetag{A.8} in terms of the coefficients $a_1$, $a_2$, $a_3$, $a_4$ 
of the initial polynomial we need to express $x_0$ through them. Let's recall 
that the equality \mythetag{A.6} means that $x_0$ is a common root of $P_4(x)$
and its first derivative $P_3(x)=P^{\kern 1pt\prime}_4(x)$:
$$
\hskip -2em
P_3(x)=4\,x^3+3\,a_1\,x^2+2\,a_2\,x+a_3. 
\mytag{A.11}
$$
Let's combine the polynomials \mythetag{A.1} and \mythetag{A.11} in the following
way: 
$$
\hskip -2em
Q_3(x)=4\,P_4(x)-x\,P_3(x). 
\mytag{A.12}
$$
It is easy to see that \mythetag{A.12} is another cubic polynomial:
$$
\hskip -2em
Q_3(x)=a_1\,x^3+2\,a_2\,x^2+3\,a_3\,x+4\,a_4. 
\mytag{A.13}
$$
Using \mythetag{A.13}, we combine it with \mythetag{A.11} as follows: 
$$
\hskip -2em
Q_2(x)=4\,Q_3(x)-a_1\,P_3(x).
\mytag{A.14}
$$
It is easy to see that \mythetag{A.14} is a quadratic polynomial:
$$
\hskip -2em
Q_2(x)=(8\,a_2-3\,a_1^2)\,x^2+(12\,a_3-2\,a_1\,a_2)\,x+16\,a_4-a_1\,a_3.
\mytag{A.15}
$$
Another quadratic polynomial is derived by means of the formula
$$
R_2(x)=(8\,a_2-3\,a_1^2)\,P_3(x)-4\,x\,Q_2(x).
\mytag{A.16}
$$
The explicit formula for the polynomial \mythetag{A.16} looks like
$$
\hskip -2em
\gathered
R_2(x)=(32\,a_1\,a_2-48\,a_3-9\,a_1^3)\,x^2+(4\,a_1\,a_3-64\,a_4\,-\\
-\,6\,a_2\,a_1^2+16\,a_2^2)\,x+8\,a_3\,a_2-3\,a_1^2\,a_3.
\endgathered
\mytag{A.17}
$$
Now we combine \mythetag{A.15} and \mythetag{A.17} by means of the formula
$$
\hskip -2em
P_1(x)=\frac{8\,a_2-3\,a_1^2}{16}\,R_2(x)-\frac{32\,a_1\,a_2-48\,a_3-9\,a_1^3}{16}
\,Q_2(x).
\mytag{A.18}
$$
It is easy to see that \mythetag{A.18} is a linear polynomial:
$$
\hskip -2em
\gathered
P_1(x)=(8\,a_2^3+36\,a_3^2+6\,a_1^3\,a_3-32\,a_2\,a_4-2\,a_1^2\,a_2^2\,+\\
+\,12\,a_1^2\,a_4-28\,a_1\,a_2\,a_3)\,x-3\,a_1\,a_3^2+48\,a_4\,a_3
+9\,a_1^3\,a_4\,+\\
+\,4\,a_3\,a_2^2-a_3\,a_2\,a_1^2-32\,a_1\,a_4\,a_2. 
\endgathered
\mytag{A.19}
$$\par 
Due to $P_4(x_0)=0$ and $P_3(x_0)=0$, from \mythetag{A.12}, \mythetag{A.14},
\mythetag{A.16}, and \mythetag{A.18} we derive that $x_0$ is a root of the linear
polynomial \mythetag{A.19}:
$$
P_1(x_0)=0. 
\mytag{A.20}
$$
Let's denote through $A_0$ the leading coefficient of the polynomial 
$P_1(x)$:
$$
A_0=8\,a_2^3+36\,a_3^2+6\,a_1^3\,a_3-32\,a_2\,a_4-2\,a_1^2\,a_2^2
+12\,a_1^2\,a_4-28\,a_1\,a_2\,a_3.
\quad
\mytag{A.21}
$$
Similarly, lets denote through $A_1$ the constant term of the polynomial 
$P_1(x)$: 
$$
A_1=-3\,a_1\,a_3^2+48\,a_4\,a_3+9\,a_1^3\,a_4+4\,a_3\,a_2^2
-a_3\,a_2\,a_1^2-32\,a_1\,a_4\,a_2.
\quad
\mytag{A.22}
$$
Due to \mythetag{A.21} and \mythetag{A.22} the equation \mythetag{A.20} is
written as 
$$
\hskip -2em
A_0\,x_0+A_1=0. 
\mytag{A.23}
$$\par 
     Let's begin with the case $A_0\neq 0$. In this case we can resolve the 
linear equation \mythetag{A.23} with respect to the variable $x_0$. Let's 
substitute 
$$
\hskip -2em
x_0=-\frac{A_1}{A_0}
\mytag{A.24}
$$
into the formula \mythetag{A.10} for the discriminant $D_2$. As a result
we get the fraction 
$$
\hskip -2em
D_2=\frac{4\,B_2}{A_0^{2}},
\mytag{A.25}
$$
where the numerator term $B_2$ is given by the following explicit formula: 
$$
\hskip -2em
\gathered
B_2=552\,a_2^2\,a_1^4\,a_3^2-30\,a_1^6\,a_4\,a_2^2-64\,a_2^7
+2208\,a_1\,a_3^3\,a_2^2\,-\\
-\,616\,a_2^3\,a_1^2\,a_3^2-704\,a_2^4\,a_1^2\,a_4
+264\,a_2^3\,a_1^4\,a_4+1536\,a_4\,a_3^2\,a_2^2\,-\\
-\,336\,a_1^3\,a_2^4\,a_3+480\,a_1\,a_2^5\,a_3+78\,a_1^5\,a_2^3\,a_3
-900\,a_2\,a_1^3\,a_3^3\,+\\
+\,144\,a_2\,a_1^4\,a_4^2-126\,a_2\,a_1^6\,a_3^2+900\,a_1^4\,a_4\,a_3^2
-1152\,a_1^3\,a_4^2\,a_3\,+\\
+\,2304\,a_1\,a_3^3\,a_4-1296\,a_2\,a_3^4-1024\,a_2^3\,a_4^2
+512\,a_2^5\,a_4\,-\\
-\,608\,a_2^4\,a_3^2-12\,a_1^4\,a_2^5+48\,a_1^2\,a_2^6
-18\,a_1^6\,a_4^2+198\,a_1^2\,a_3^4\,-\\
-\,4608\,a_4^2\,a_3^2+90\,a_1^5\,a_3^3+a_1^6\,a_2^4+9\,a_1^8
\,a_3^2+2112\,a_1^3\,a_3\,a_2^2\,a_4\,-\\
-\,1024\,a_1\,a_3\,a_2^3\,a_4-4032\,a_2\,a_1^2\,a_3^2\,a_4
-828\,a_2\,a_1^5\,a_4\,a_3\,+\\
+\,4608\,a_4^2\,a_3\,a_1\,a_2+90\,a_1^7\,a_4\,a_3-6\,a_1^7\,a_2^2\,a_3.
\endgathered
\mytag{A.26}
$$
Since the denominator in \mythetag{A.25} is positive, the inequality 
\mythetag{A.8} can be written as 
$$
\hskip -2em
B_2\leqslant 0. 
\mytag{A.27} 
$$
Summarizing the above calculations, now we can formulate a lemma.
\mylemma{A.1} If a quartic polynomial $P_4(x)$ with real coefficients 
in \mythetag{A.1} has exactly one real root and if $A_0\neq 0$ in 
\mythetag{A.21}, then the discriminant $D_4=0$ in \mythetag{A.7} and 
$B_2\leqslant 0$ in \mythetag{A.26}. 
\endproclaim
The case $B_2=0$ in \mythetag{A.27} or, equivalently, the case $D_2=0$ 
in \mythetag{A.8} is exceptional. Omitting this case, we can formulate 
a lemma converse to Lemma~\mythelemma{A.1}. 
\mylemma{A.2} If $D_4=0$ in \mythetag{A.7}, if $A_0\neq 0$ in \mythetag{A.21}
and if $B_2<0$ in \mythetag{A.26}, then the quartic polynomial $P_4(x)$ 
with real coefficients in \mythetag{A.1} has exactly one real root of 
multiplicity $2$. 
\endproclaim
     Let's proceed to the case $A_0=0$ assuming the condition $D_4=0$ 
in \mythetag{A.6} is fulfilled. The condition $D_4=0$ means that the polynomial
$P_4(x)$ has a root $x_0$ of multiplicity at least two. Hence we can apply 
\mythetag{A.2} and the formulas \mythetag{A.4} and \mythetag{A.5}
following from \mythetag{A.2} and \mythetag{A.3}. The formulas 
\mythetag{A.4} can be written as 
$$
\xalignat 2
&\hskip -2em
a_1=b_1-2\,x_0,
&&a_2=b_2+x_0^2-2\,b_1\,x_0.
\mytag{A.28}
\endxalignat
$$
Substituting \mythetag{A.28} and \mythetag{A.5} into \mythetag{A.21}, we derive 
$$
\hskip -2em
A_0=-2\,(b_1^2-4\,b_2)\,(x_0^2+b_1\,x_0+b_2)^2.
\mytag{A.29}
$$
Comparing the formula \mythetag{A.29} with the formulas \mythetag{A.3} and 
\mythetag{A.9}, we see that the equality $A_0=0$ implies at least one of the 
following two equalities:
$$
\xalignat 2
&\hskip -2em
D_2=0,
&&P_2(x_0)=0.
\mytag{A.30}
\endxalignat
$$
The equality $D_2=0$ means that the quadratic polynomial $P_2(x)$ has a double
root: 
$$
\hskip -2em
P_2(x)=(x-x_1)^2. 
\mytag{A.31}
$$
The equality $P_2(x_0)=0$ means that $x_0$ is one of two roots of the 
quadratic polynomial $P_2(x)$, i\.\,e\. the polynomial $P_2(x)$ is written as
$$
\hskip -2em
P_2(x)=(x-x_0)\,(x-x_1). 
\mytag{A.32}
$$
In both cases \mythetag{A.31} and \mythetag{A.32} the quartic polynomial
$P_4(x)=(x-x_0)^2\,P_2(x)$ has exactly one real root if and only if $x_1=x_0$. 
Applying $x_1=x_0$ back to \mythetag{A.31} and \mythetag{A.32}, we derive 
$D_2=0$ and $P_2(x_0)=0$, i\.\,e\. both equalities \mythetag{A.30} are fulfilled.
In other words we have the following lemma. 
\mylemma{A.3} If $A_0=0$ in \mythetag{A.21}, then the equality $D_4=0$ in 
\mythetag{A.7} and the couple of equalities \mythetag{A.30} constitute a necessary 
and sufficient condition for the quartic polynomial $P_4(x)$ with real coefficients 
in \mythetag{A.1} to have exactly one real root, which is of multiplicity $4$ in
this case. 
\endproclaim
     Note that the equalities \mythetag{A.30} are not written in terms of the 
coefficients $a_1$, $a_2$, $a_3$, $a_4$ of the polynomial $P_4(x)$. In order to 
write them properly we apply $x_1=x_0$ to the equalities \mythetag{A.31} and 
\mythetag{A.32} and derive from them 
$$
\xalignat 2
&\hskip -2em
P_2(x)=(x-x_0)^2,
&&P_4(x)=(x-x_0)^4.
\mytag{A.33}
\endxalignat
$$
The equalities \mythetag{A.33} mean that the polynomials $P_2(x)$ and $x-x_0$ are 
expressed through the derivatives of the polynomial $P_4(x)$:
$$
\xalignat 2
&\hskip -2em
P_2(x)=\frac{P_4^{\kern 1pt\prime\prime}(x)}{12},
&&x-x_0=\frac{P_4^{\kern 1pt\prime\prime\prime}(x)}{24}.
\mytag{A.34}
\endxalignat
$$
Substituting \mythetag{A.1} into \mythetag{A.34}, we derive 
$$
\xalignat 2
&\hskip -2em
P_2(x)=x^2+\frac{a_1}{2}\,x+\frac{a_2}{6}
&&x-x_0=x+\frac{a_1}{4}.
\mytag{A.35}
\endxalignat
$$
Using \mythetag{A.3}, from \mythetag{A.35} we derive the formulas
for $b_1$, $b_2$, and $x_0$:
$$
\xalignat 3
&\hskip -2em
b_1=\frac{a_1}{2},
&&b_2=\frac{a_2}{6},
&&x_0=-\frac{a_1}{4}.
\quad
\mytag{A.36}
\endxalignat
$$
The third formula  \mythetag{A.36} replaces the formula \mythetag{A.24}
in the present case $A_0=0$. Substituting it into \mythetag{A.10}, we 
derive the formula for $D_2$: 
$$
\hskip -2em
D_2=\frac{3}{2}\,a_1^2-4\,a_2.
\mytag{A.37}
$$
The formula  \mythetag{A.37} replaces the formula \mythetag{A.25}
in the present case $A_0=0$.\par
     In order to express $P_2(x_0)$ through $a_1$, $a_2$, $a_3$, $a_4$
we first apply \mythetag{A.4} to \mythetag{A.3}. As a result we get
the following expression for $P_2(x)$:
$$
\hskip -2em
P_2(x)=x^2+(a_1+2\,x_0)\,x+a_2+3\,x_0^2+2\,x_0\,a_1.
\mytag{A.38}
$$
Then we substitute $x=x_0$ into \mythetag{A.38} and get 
$$
\hskip -2em
P_2(x_0)=6\,x_0^2+3\,x_0\,a_1+a_2.
\mytag{A.39}
$$
And finally we apply the third equality \mythetag{A.36} to \mythetag{A.39}.
This yields
$$
\hskip -2em
P_2(x_0)=-\frac{3}{8}\,a_1^2+a_2.
\mytag{A.40}
$$
Comparing \mythetag{A.40} with \mythetag{A.37}, we see that two equalities 
\mythetag{A.30} become equivalent to each other. Lemma~\mythelemma{A.3}
now is reformulated as follows. 
\mylemma{A.4} If $A_0=0$ in \mythetag{A.21}, then the equality $D_4=0$ in 
\mythetag{A.7} and the equality $D_2=0$ in \mythetag{A.37} constitute a necessary 
and sufficient condition for the quartic polynomial $P_4(x)$ with real coefficients 
in \mythetag{A.1} to have exactly one real root, which is of multiplicity $4$ in
this case. 
\endproclaim
\mytheorem{A.1} A quartic polynomial $P_4(x)$ with real coefficients 
in \mythetag{A.1} has exactly one real root if and only if its discriminant
$D_4=0$ in \mythetag{A.7} and if one of the following two conditions is fulfilled:
\roster
\item"1)" $A_0\neq 0$ in \mythetag{A.21} and $B_2<0$ in \mythetag{A.26};
\item"2)" $A_0=0$ in \mythetag{A.21} and $D_2=0$ in \mythetag{A.37}.
\endroster
\endproclaim
     Theorem~\mythetheorem{A.1} is the ultimate result. It summarizes 
Lemmas~\mythelemma{A.1}, \mythelemma{A.2}, and \mythelemma{A.4}.
and provides a necessary and sufficient condition for a real quartic
polynomial to have exactly one real root.\par

\enddocument
\end